\documentclass[english]{article}
\usepackage[T1]{fontenc}
\usepackage[latin9]{inputenc}
\usepackage{float}
\usepackage{textcomp}
\usepackage{amsmath}
\usepackage{amssymb}
\usepackage{graphicx}
\usepackage{babel}
\begin{document}

\title{Non-local image deconvolution by Cauchy sequence}

\author{Jack Dyson, Gianni Albertini}
\maketitle
\begin{abstract}
Abstract: We present the deconvolution between two smooth function
vectors as a Cauchy sequence of weight functions. From this we develop
a Taylor series expansion of the convolution problem that leads to
a non-local approximation for the deconvolution in terms of continuous
function spaces. Optimisation of this form against a given measure
of error produces a theoretically more exact algorithm. The discretization
of this formulation provides a deconvolution iteration that deconvolves
images quicker than the Richardson-Lucy algorithm.
\end{abstract}
\textbf{Keywords}: deconvolution, image processing, Banach spaces,
Cauchy sequences, Richardson-Lucy

\section*{Introduction}

The deconvolution of images is a frequent inverse problem occurring
in scientific and engineering applications. In particular, problems
involving smooth images and kernels in two or more dimensions are
described by:

\begin{equation}
h=\int k(\boldsymbol{x}-\boldsymbol{x}')g(\boldsymbol{x}')d\boldsymbol{x}'\label{eq:1}
\end{equation}

The objective, given $h$ and $k$, is to find the image $g$ where
$h$ is the observed image, $g$ is the function to be determined
and $k$ is the convolving kernel. We explicitly assume in the following
analysis that the images are differentiable and noiseless.

Given two continuous integrable images, $h$ and $g$ that are related
by equation (1), it can be shown that there exists a continuous bijection
on the complete Banach space of functions $\mathcal{L}_{2}$ that
will deconvolve $g$ from $h$. Writing:

\begin{equation}
h=g-(\mathcal{I}-k*)g\label{eq:2}
\end{equation}

the solution $g$ can be written:

\begin{equation}
g=\sum_{n=0}^{\infty}\left(\mathcal{I}-k*\right)^{n}h\label{eq:3}
\end{equation}

Equation (\ref{eq:3}) is the Von Neumann series from the general
theory of linear operators \cite{key-1,key-3}. However, the convolution
integral has \emph{additional} mathematical structure that can be
used to optimise the deconvolution beyond the general iterative process
given in equation (\ref{eq:3}). The approach of this article is to
determine an alternative convergent series for equation (\ref{eq:1}).
Section \ref{sec:1} presents an alternative to equation (\ref{eq:3}).
Section \ref{sec:2} develops the series given in section \ref{sec:1}.
In section \ref{sec:3} the iteration developed in section \ref{sec:2}
is accelerated using an auxilliary Richardson-Lucy iteration. Section
\ref{sec:4} develops a simple noise suppression strategy for the
results of section \ref{sec:3}. Section \ref{sec:5} numerically
analyses these results by directly deconvolving a discrete test image.
Section \ref{sec:6} concludes. Section \ref{sec:7} presents figures.
Section \ref{sec:8} cites the most important references used in the
development of this work.

\section{Existence of an alternative parameterisation\label{sec:1}}

From equation (\ref{eq:3}) derive a sequence of partial sums $(g_{n}):g_{n}\in\mathcal{L}_{2}$.
Define the sequence of weight functions $(\varrho_{n}):g_{n}h^{-1},h\ne0\implies$from
equation (\ref{eq:1}):$\|k*\|_{2}^{-1}\le\|\varrho_{n}\|_{2}\sim1$.
Assuming that $\|\mathcal{I}-k*\|_{2}<1$ in equation (\ref{eq:3})
$\implies(g_{n})$ and $(\varrho_{n})$ are Cauchy sequences in the
Banach space. Note in particular that for problems where $k*$ is
close to the identity operation, the $n$th weight $\varrho_{n}\sim1$.
Next we apply an important result to the sequence $(g_{n})$:

\textbf{Lemma 1\label{LS1.1}} \textit{The convergent sequence $(g_{n})$
of sums is not unique.} Proof: The sequence $(g_{n})$ is by definition
absolutely convergent. Therefore the series of some subsequence of
$a_{n}=\sup_{m>n}g_{n}-g_{m}$ is also convergent and rearrangable.
Writing $g_{n_{i}}=\sum_{k=0}^{i}a_{n_{k}}+g_{n_{i+1}}$ generates
an arbitrary subsequence of $(g_{n})$ which proves the lemma. 

Using $g=\varrho h$ in equation (\ref{eq:1})$\implies h=\int k(\boldsymbol{x}-\boldsymbol{x}')\varrho(\boldsymbol{x}')h(\boldsymbol{x}')d\boldsymbol{x}'$.
Hence, if $k*\rightarrow\mathcal{I},\varrho\rightarrow1\therefore g\rightarrow h^{2}\left[k*h\right]^{-1}$,
and substituting in equation (\ref{eq:3}) obtains:
\begin{equation}
\sum_{n=1}^{\infty}\left(\mathcal{I}-k\right)^{n}h\rightarrow\left(1-\frac{h}{k*h}\right)h\label{eq:4}
\end{equation}

\pagebreak{}Combining Lemma \ref{LS1.1}, equations (\ref{eq:4})
and (\ref{eq:3}) the following is true:

\begin{equation}
\varrho=\varrho_{0}+\sum_{n}t_{n}\label{eq:5}
\end{equation}

where:

\begin{equation}
\varrho_{0}=\frac{h}{k*h}\label{eq:5-1}
\end{equation}

The sequence $\left(\varrho_{n}\right)$ where $\varrho_{n}=\varrho_{0}+\sum_{k=1}^{n}t_{k}$
is an alternative to the Von Neumann series for $(g_{n})$, where
both $\lim h\varrho_{n}$ and $\lim g_{n}$ are equal to the solution
$g$. In the next section we derive an explicit form for equation
(\ref{eq:5}).

\section{Explicit form for an alternative sequence\label{sec:2}}

Given equation (\ref{eq:5}) note the well known property of convolution
integrals, $k*h=h*k$. Together with $g=\varrho h$ we have $h=\int k\left(\boldsymbol{x}'\right)\varrho(\boldsymbol{x}-\boldsymbol{x}')h(\boldsymbol{x}-\boldsymbol{x}')d\boldsymbol{x}'$.
Expanding $\varrho(\boldsymbol{x}-\boldsymbol{x}')$ \emph{about}
$\boldsymbol{x}$:

\begin{equation}
\varrho(\boldsymbol{x}-\boldsymbol{x}')=\sum_{n=0}^{\infty}\frac{1}{n!}(-\boldsymbol{x}'\cdot\nabla_{\boldsymbol{x}})^{n}\varrho(\boldsymbol{x})\label{eq:11}
\end{equation}

Used in the convolution integral, equation (\ref{eq:11}) factors
out the derivatives of $\varrho$ \cite{key-8,key-11}. 

\begin{equation}
h=\left(\hat{\boldsymbol{m}}\cdot\nabla_{\boldsymbol{x}}\right)^{0}\varrho+\left(\hat{\boldsymbol{m}}\cdot\nabla_{\boldsymbol{x}}\right)^{1}\varrho+\frac{1}{2!}\left(\hat{\boldsymbol{m}}\cdot\nabla_{\boldsymbol{x}}\right)^{2}\varrho+\cdots\label{eq:12}
\end{equation}

where we define the following operator symbols for $\boldsymbol{x}=(x_{1},x_{2})$:

\[
\left(\hat{\boldsymbol{m}}\cdot\nabla_{\boldsymbol{x}}\right)^{0}\equiv k*h
\]

\[
\hat{m}_{d}^{n}\left(\boldsymbol{x}\right)\equiv(-1)^{n}\int(x_{d}-x_{d}')^{n}k\left(\boldsymbol{x}-\boldsymbol{x}'\right)h(\boldsymbol{x}')d\boldsymbol{x}'
\]

\[
\hat{m}_{1}^{n-k}\hat{m}_{2}^{k}\left(\boldsymbol{x}\right)\equiv(-1)^{n}\int(x_{1}-x_{1}')^{n-k}\left(x_{1}-x_{2}'\right)^{k}k\left(\boldsymbol{x}-\boldsymbol{x}'\right)h(\boldsymbol{x}')d\boldsymbol{x}'
\]

so that $\left(\hat{\boldsymbol{m}}\cdot\nabla_{\boldsymbol{x}}\right)^{n}=\sum_{k}C_{k}^{n}\hat{m}_{1}^{n-k}\hat{m}_{2}^{k}\left(\boldsymbol{x}\right)\partial_{1}^{n-k}\partial_{2}^{k}$,
$d=1,2$ and $k=0,1,2,\dots n$. The scalar product operator $\hat{\boldsymbol{m}}\left(\boldsymbol{x}\right)\cdot=\left(\hat{m}_{1}\left(\boldsymbol{x}\right),\hat{m}_{2}\left(\boldsymbol{x}\right)\right)\cdot$
is represented as a vector integral operation on the nabla, $\nabla_{\boldsymbol{x}}$.
Noting that the antisymmetric moments \cite{key-8} $\hat{m}_{1}^{n-k}\hat{m}_{2}^{k}\left(\boldsymbol{x}\right)\approx0$
for $n,k\notin\left\{ 2,4,\dots\right\} $ together with Stirling's
approximation for $2n!$ leads to a direct representation of the convolution
in terms of the weight functions $\varrho_{0}$ and $\varrho$:

\begin{equation}
\varrho_{0}=\sum_{n=0}^{\infty}\frac{1}{n!}\left[\hat{\boldsymbol{M}}_{n}\cdot\nabla_{\boldsymbol{x}}\right]^{2n}\varrho\leq\sum_{n=0}^{\infty}\frac{1}{n!}\left[\hat{\boldsymbol{m}}\cdot\nabla_{\boldsymbol{x}}\right]^{2n}\varrho=e^{\left(\hat{\boldsymbol{m}}\cdot\nabla_{\boldsymbol{x}}\right)^{2}}\varrho\label{eq:13}
\end{equation}

where $\hat{\boldsymbol{M}}_{n}\left(\boldsymbol{x}\right)\propto n^{-1}\left[k*h\right]^{-n^{-1}}\hat{\boldsymbol{m}}\left(\boldsymbol{x}\right)$.
The inequality becomes an equality in the limit $\left\Vert \hat{\boldsymbol{m}}\left(\boldsymbol{x}\right)\right\Vert _{2}\rightarrow0$.
The integral moment ratios in equation (\ref{eq:13}) vary slowly
over $h$ and in particular for symmetric kernels we have $\hat{m}_{1}^{2}\left(\boldsymbol{x}\right)=\hat{m}_{2}^{2}\left(\boldsymbol{x}\right)=\hat{m}^{2}\left(\boldsymbol{x}\right)$.
Using these facts in equation (\ref{eq:13}) gives :

\[
\varrho_{0}\leq e^{\alpha^{2}\nabla_{\boldsymbol{x}}^{2}}\varrho,\ \alpha\in\mathbb{R}
\]

If the terms of the rightmost sum in equation (\ref{eq:13}) are negligible
beyond $n=1$ then:

\begin{equation}
\varrho_{0}=e^{\alpha^{2}\nabla_{\boldsymbol{x}}^{2}}\varrho\label{eq:14}
\end{equation}

represents the convolution well for a fixed $\alpha$ value. Inverting
this relationship estimates the deconvolution :

\begin{equation}
\varrho=e^{-\alpha^{2}\nabla_{\boldsymbol{x}}^{2}}\varrho_{0}\label{eq:15}
\end{equation}

By definition:

\begin{equation}
e^{-\alpha^{2}\nabla_{\boldsymbol{x}}^{2}}=\lim_{n\rightarrow\infty}\left(1-\frac{\alpha^{2}\nabla_{\boldsymbol{x}}^{2}}{n}\right)^{n}\label{eq:16}
\end{equation}
which implies the approximate iteration formula:

\begin{equation}
\varrho_{n}=\varrho_{n-1}-\frac{\alpha^{2}}{n}\nabla_{\boldsymbol{x}}^{2}\varrho_{n-1}\label{eq:17}
\end{equation}

The inverted convolution sequence $\left(\varrho_{n}\right)$ is bounded
by the exponential series in equation (\ref{eq:13}) and is therefore
Cauchy, since equation (\ref{eq:16}) $\implies\left\Vert \varrho_{n}-\varrho_{n-1}\right\Vert _{2}\rightarrow0$.
Further, since equation (\ref{eq:17}) considers all the moments in
equation (\ref{eq:13}) to be equal, it appears to be \emph{local}
in the sense of \cite{key-8}: note however that the division of $\alpha$
by $n$ in equation (\ref{eq:17}) negates such a comparision.

\pagebreak{}The leftmost sum of equation (\ref{eq:13}) would be
better reproduced by iteration (\ref{eq:17}) when $\alpha$ is considered
to be an unknown function of $n$, since iteration (\ref{eq:17})
is an upper bound approximation to the correct sequence $\left(\varrho_{n}\right)$.
Therefore, given that the convolution problem is \emph{not} specified
by a unique Cauchy sequence (lemma \ref{LS1.1}), equation (\ref{eq:17})
can be made \emph{optimal} relative to some measure of error if we
write:

\begin{equation}
\varrho_{n}=\varrho_{n-1}-\frac{\alpha_{n}^{2}}{n}\nabla_{\boldsymbol{x}}^{2}\varrho_{n-1}\label{eq:20}
\end{equation}

for a suitable choice of the infinite set $\left\{ \alpha_{n}\right\} $
- see equation (\ref{eq:32}). Equation (\ref{eq:20}) therefore solves
the deconvolution problem in the sense of equation (\ref{eq:5}).

\section{Auxilliary acceleration\label{sec:3}}

For large $n$, to a good approximation, we have $h\varrho_{n+1}=h\varrho_{n}\implies$
substitution for $h\approx k*\left(h\varrho_{n}\right)$ on the left
hand side of this expression leads to the well known Richardson-Lucy
\cite{key-4} propagator for $\varrho$:

\begin{equation}
\varrho_{n+1}=\frac{h\varrho_{n}}{k*\left(h\varrho_{n}\right)}\label{eq:22}
\end{equation}

Combining equations (\ref{eq:22}) and (\ref{eq:20}) where the correction
in the denominator is negligible for large $n$:

\begin{equation}
\varrho_{n+1}=\frac{h}{k*\left(h\varrho_{n-1}\right)}\left[\varrho_{n-1}-\frac{\alpha_{n}^{2}}{n}\nabla_{\boldsymbol{x}}^{2}\varrho_{n-1}\right]\label{eq:25}
\end{equation}

\section{Noise suppression\label{sec:4}}

Equation (\ref{eq:22}) is recognisable as the $n$th iteration of
the Richardson-Lucy algorithm without the accompanying right convolution
with $k$. Therefore it is prudent to modify equation (\ref{eq:22}):

\begin{equation}
\varrho_{n+1}=\left[\frac{h}{k*\left(h\varrho_{n-1}\right)}*k\right]\varrho_{n}\label{eq:26}
\end{equation}

Using this in equation (\ref{eq:20}):

\begin{equation}
\varrho_{n+1}=\frac{h}{k*\left(h\varrho_{n-1}\right)}*k\left[\varrho_{n-1}-\frac{\alpha_{n}^{2}}{n}\nabla_{\boldsymbol{x}}^{2}\varrho_{n-1}\right]\label{eq:27}
\end{equation}

Starting the algorithm with an initial value $\varrho_{0}$ will generate
a nabla term acting on every $\varrho_{n}$ that generates noise in
the solution at every iteration. This build up of noise is proportional
to $n$ and destroys the convergence properties of the algorithm beyond
a certain number of iterations.

The noise statistics of a differential operator exhibit a \textquotedblleft random
walk\textquotedblright{} phenomenon in its function space \emph{without}
any overall cancellation in the noise. After a sufficient number of
iterations the norm of the derivative is characterised just by the
noise function $\phi_{n}$ and, $\left\Vert \nabla_{\boldsymbol{x}}^{2}\right\Vert _{2}=\left\Vert \phi_{n}\right\Vert _{2}\propto n^{2}$.
Therefore iteration (\ref{eq:27}) is no longer representative of
the exponential operator series (\ref{eq:17}) when noise is present,
rather the operator diverges as $n$ for large values of $n$. It
is therefore prudent to normalise the derivative $\nabla_{\boldsymbol{x}}^{2}\rightarrow\nabla_{\boldsymbol{x}}^{2}/\left\Vert \nabla_{\boldsymbol{x}}^{2}\right\Vert _{2}$.

This strategy gives good noise suppression as intended, which suggests
that the power of $n$ in the denominator can be tuned just enough
to correct for derivative noise \emph{and} at the same time optimise
the signal propagator so that overall the algorithm converges correctly
according to equation (\ref{eq:20}). Using these developments in
equation (\ref{eq:27}) gives the following iteration for a 2 dimensional
image convolved by a symmetric kernel $k$:

\begin{equation}
\varrho_{n+1}=\frac{h}{k*\left(h\varrho_{n}\right)}*k\left[\varrho_{n-1}-\frac{\alpha_{n}^{2}}{\left\Vert \nabla_{\boldsymbol{x}}^{2}\right\Vert _{2}^{p}}\nabla_{\boldsymbol{x}}^{2}\varrho_{n-1}\right]\label{eq:31}
\end{equation}

where $p\approx1$ and $\alpha_{n}$ need to be determined per iteration
in order that convergence is optimised - see equation (\ref{eq:32}).

\section{Results and observations\label{sec:5}}

All results are computed with MATLAB\textsuperscript{©} and the algorithm
is implemented in discrete form by modifying the MATLAB\textsuperscript{©}
``deconvlucy'' command in accordance with equation (\ref{eq:31})
multiplied through by $h$. All norms are taken as root mean square
(RMS) values.

\begin{figure}[H]
\noindent \begin{centering}
\includegraphics[width=0.333333\textwidth,height=0.333333\textwidth]{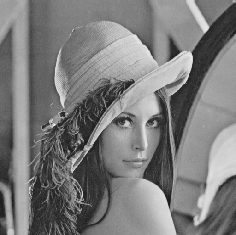}\includegraphics[width=0.333333\textwidth,height=0.33333\textwidth]{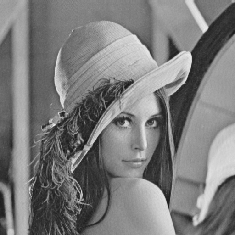}\includegraphics[width=0.333333\textwidth,height=0.333333\textwidth]{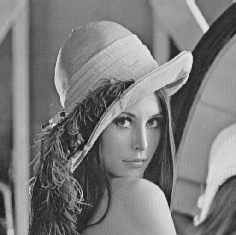}
\par\end{centering}

\noindent \begin{centering}
\includegraphics[clip,width=0.333333\textwidth,height=0.333333\textwidth]{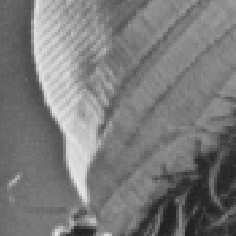}\includegraphics[width=0.333333\textwidth,height=0.333333\textwidth]{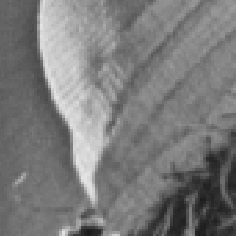}\includegraphics[width=0.333333\textwidth,height=0.333333\textwidth]{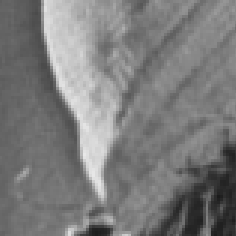}
\par\end{centering}

\caption{(Left column) The \textquotedblleft Lena\textquotedblright{} original
image. (Middle column) The result of using equation \ref{eq:31} after
256 revolutions. The image indicates a region of detailed (high frequency)
restoration. (Right column) The Richardson-Lucy result after 256 revolutions.
The image blowup on the top row shows a region of less detailed restoration
compared to equation (\ref{eq:31}).}

\label{fig:1}
\end{figure}

\begin{figure}[H]
\noindent \begin{centering}
\includegraphics[width=0.5\textwidth,height=0.5\textwidth]{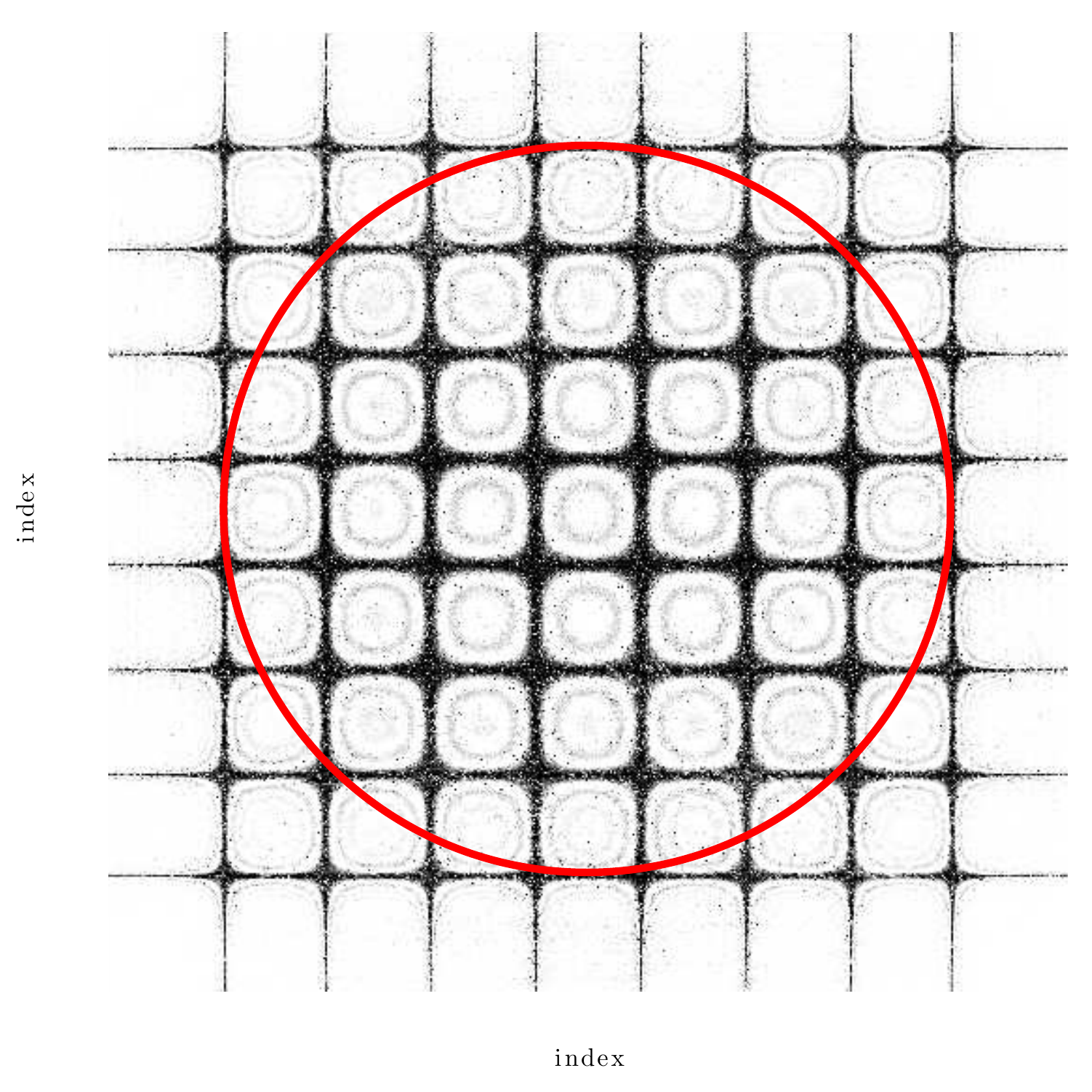}\includegraphics[width=0.5\textwidth,height=0.5\textwidth]{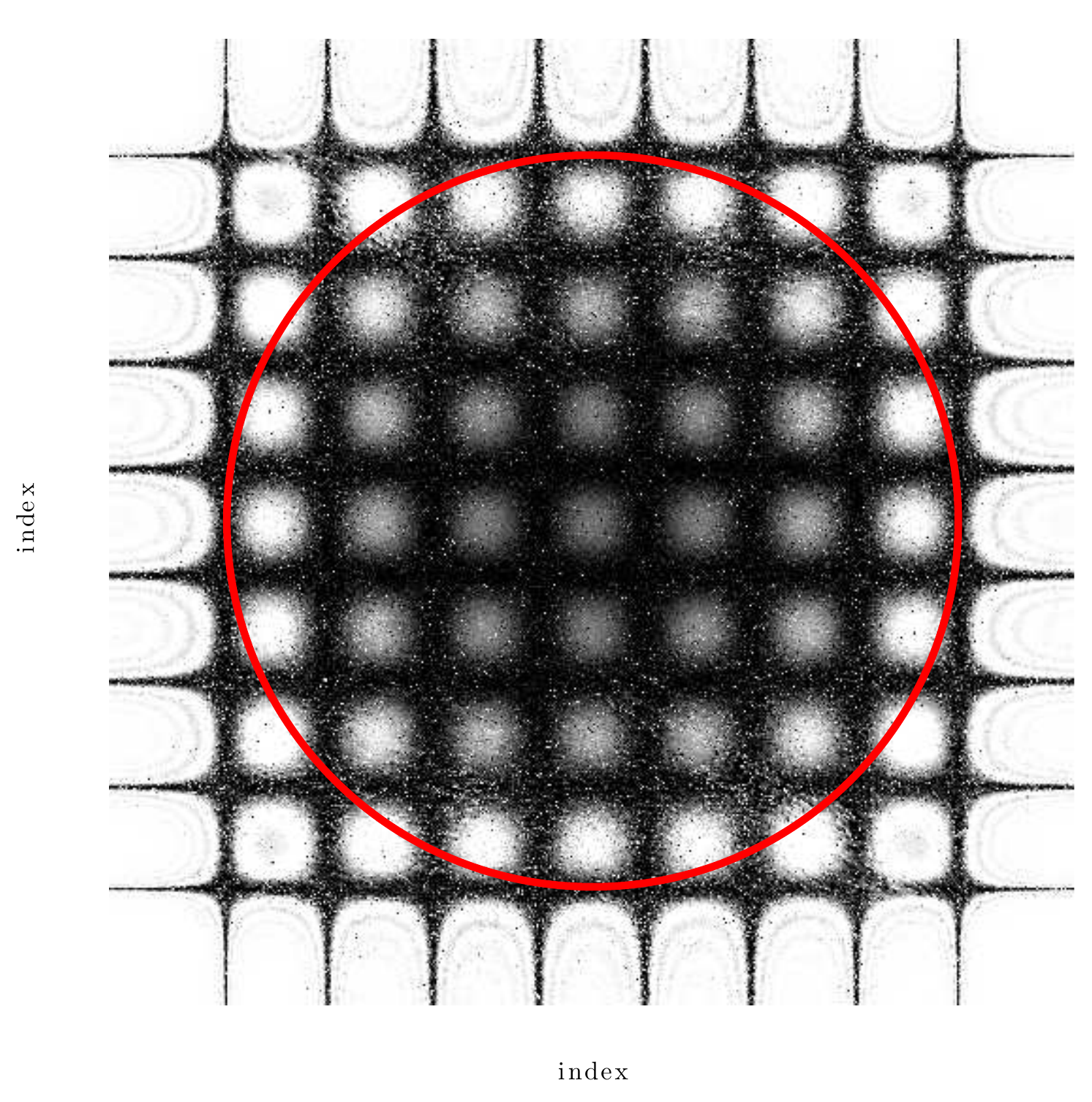}
\par\end{centering}

\caption{(Left) The final FTR image for equation (\ref{eq:31}) at 256 iterations.
The circle encompasses the high frequency regions towards the center
of the image. (Right) The final FTR image for RL at 256 iterations.}

\label{fig:3}
\end{figure}

To demonstrate the approach we convolve a $512\times512$ noiseless,
greyscale, \textquotedblleft Lena\textquotedblright{} image matrix,
$g$, by a Gaussian radial kernel of standard deviation 5 pixels clipped
at 9 pixels and then reverse the process by using equation (\ref{eq:31}).
We also take the same blurred image and deblur using the original
MATLAB\textsuperscript{©} ``deconvlucy'' command and compare the
obtained relative error in the two approaches. The results are given
in figures \ref{fig:1} and \ref{fig:3}. These show the original
image followed by the approximations obtained after 256 iterations
of each algorithm. In each case visual convergence is good. However
upon closer inspection it can be seen that equation (\ref{eq:31})
recovers the ridge detail of Lena's hat per iteration better than
the Richardson-Lucy algorithm.

Analysing the relative error, $\left\Vert \epsilon_{n}\right\Vert _{2}=\left\Vert g-g_{n}\right\Vert _{2}\left\Vert g\right\Vert _{2}^{-1}$
, it is clear that for this kind of image the present algorithm performs
better. The form of iteration (\ref{eq:20}) guarantees this because
it, when inverted, tends to reproduce the \emph{entire} leftmost series
expansion (\ref{eq:13}) for $\varrho$.

Regularisation of the measured image gives an estimate for the relative
error function:

\begin{equation}
\epsilon_{n}\left(\boldsymbol{x}\right)=1-k*\left[\varrho_{n}\left(\alpha_{1},\alpha_{2},\dots\alpha_{n}\right)h\right]h^{-1}\label{eq:32}
\end{equation}

The norm of this expression can be minimised at every iteration, to
generate an optimal $\left\{ \alpha_{n}\right\} $ set for equation
(\ref{eq:20}). Equation (\ref{eq:20}) is then made \emph{optimal}
to the $n$th order with respect to equation (\ref{eq:13}) in the
sense of the Cauchy sequence $\epsilon_{0},\epsilon_{1},\dots,\epsilon$
where $\epsilon$ is the minimum error.

We remark that the running optimisation procedure for $\alpha$ is
still the subject of ongoing research. The results presented in this
paper are for \emph{constant} $\alpha$ values.

Taking Fourier transforms of the approximate series (\ref{eq:15})
\cite{key-5} leads to an expression for the Fourier transform $\hat{\varrho}$
in terms of the Fourier transform of $\hat{\varrho}_{0}$:

\begin{equation}
\hat{\varrho}=e^{-\alpha^{2}\boldsymbol{\omega}^{2}}\hat{\varrho}_{0}\label{eq:33}
\end{equation}

which demonstrates that the operator series (\ref{eq:15}) effectively
regularises the sequence starting point $\varrho_{0}$ by an exponential
upper bound. Referring to iteration (\ref{eq:17}), the amount of
exponential regularisation applied to the Cauchy sequence is directly
controlled by the convolution moments of the image in frequency space.
Equation (\ref{eq:33}) shows that the impact on the overall ``reach''
of the algorithm into frequency space of the weight functions is characterised
by:

\begin{equation}
\omega_{R}=1/\alpha\label{eq:34}
\end{equation}

By itself, there is nothing new in an exponential regularisation,
it is already well known. The difference here is in how it arises
from Cauchy sequences, is used in equation (\ref{eq:31}) and what
it actually influences in the $\left(\varrho_{n}\right)$ sequence.
These features may make the approach useful from both a theoretical
as well as practical standpoint.

For this reason we also add another test measure, based on the norm
of the ratio of the Fourier transforms (FTR) of each image, $\text{FTR}=1-\left\Vert \hat{g}_{n}\right\Vert _{2}\left\Vert \hat{g}\right\Vert _{2}^{-1}$
where $1-\hat{g}_{n}/\hat{g}$ is proportional to $n^{-1}$ for large
$n$. 

From a practical point of view, this is is perhaps the most important
number since a value of the FTR $\sim1$ indicates that the deconvolution
algorithm is inefficient at recovering high frequency detail effectively
and acts like a low pass filter. For the theory of this paper, the
low pass filter interpretation is supported by the implication of
a small value for $\omega_{R}$ and consequently a large value for
$\alpha$ which would tend to invalidate equation (\ref{eq:14}).
In that situation the optimised form of the iteration would then be
required to make further progress.

If instead the FTR $\ll1$ for a recovered image then the chances
are that fine details in the image have been recovered well and the
algorithm consequently has good reach into the high frequency components
of the Fourier space. Theoretically, this implies that $\omega_{R}$
be large and the $\alpha$ value consequently small. In this situation
the exponential series (\ref{eq:14}) represents the deconvolution
well.

\begin{figure}[H]
\noindent \begin{centering}
\includegraphics[width=0.5\textwidth,height=0.5\textwidth]{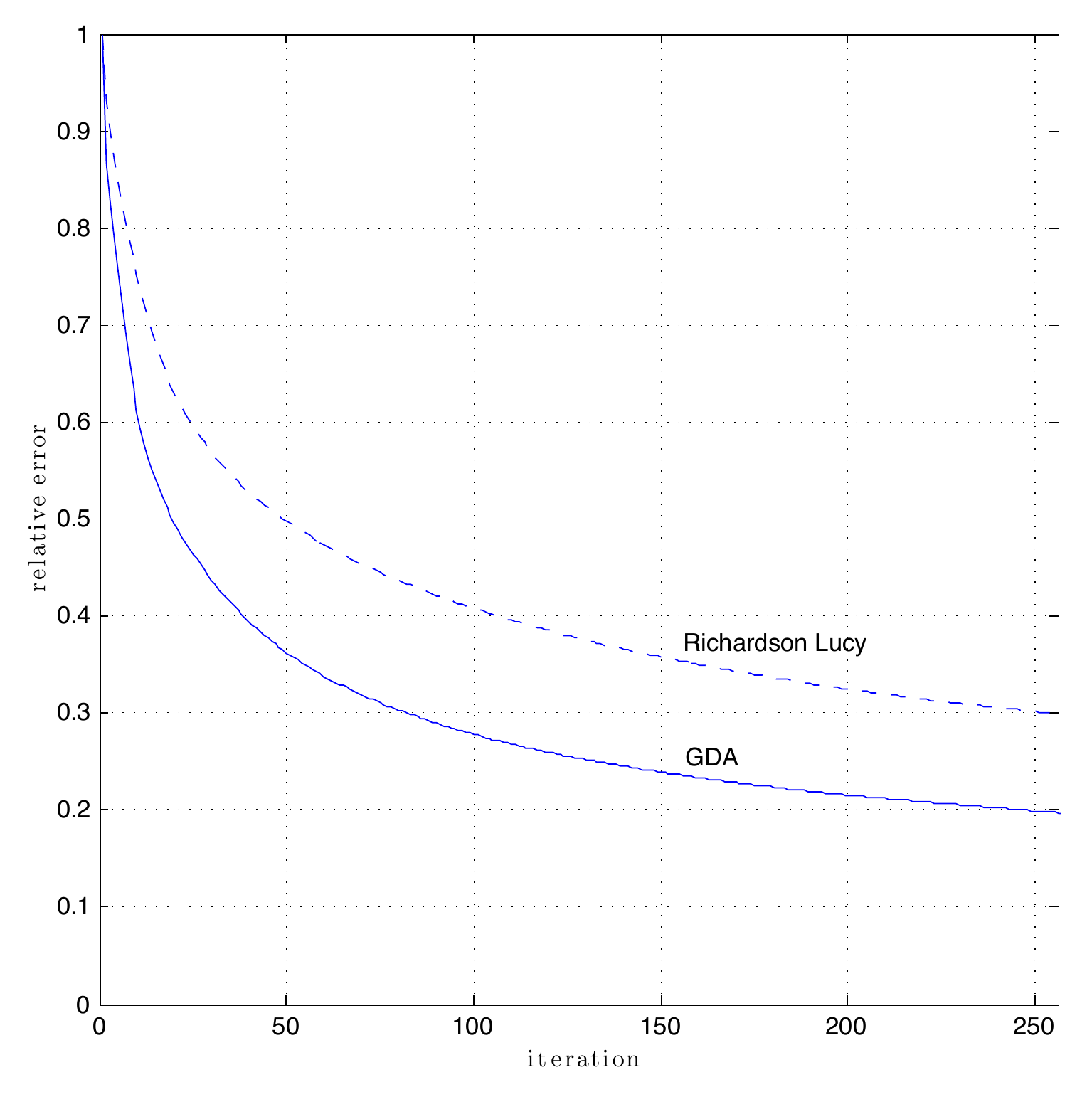}\includegraphics[width=0.5\textwidth,height=0.5\textwidth]{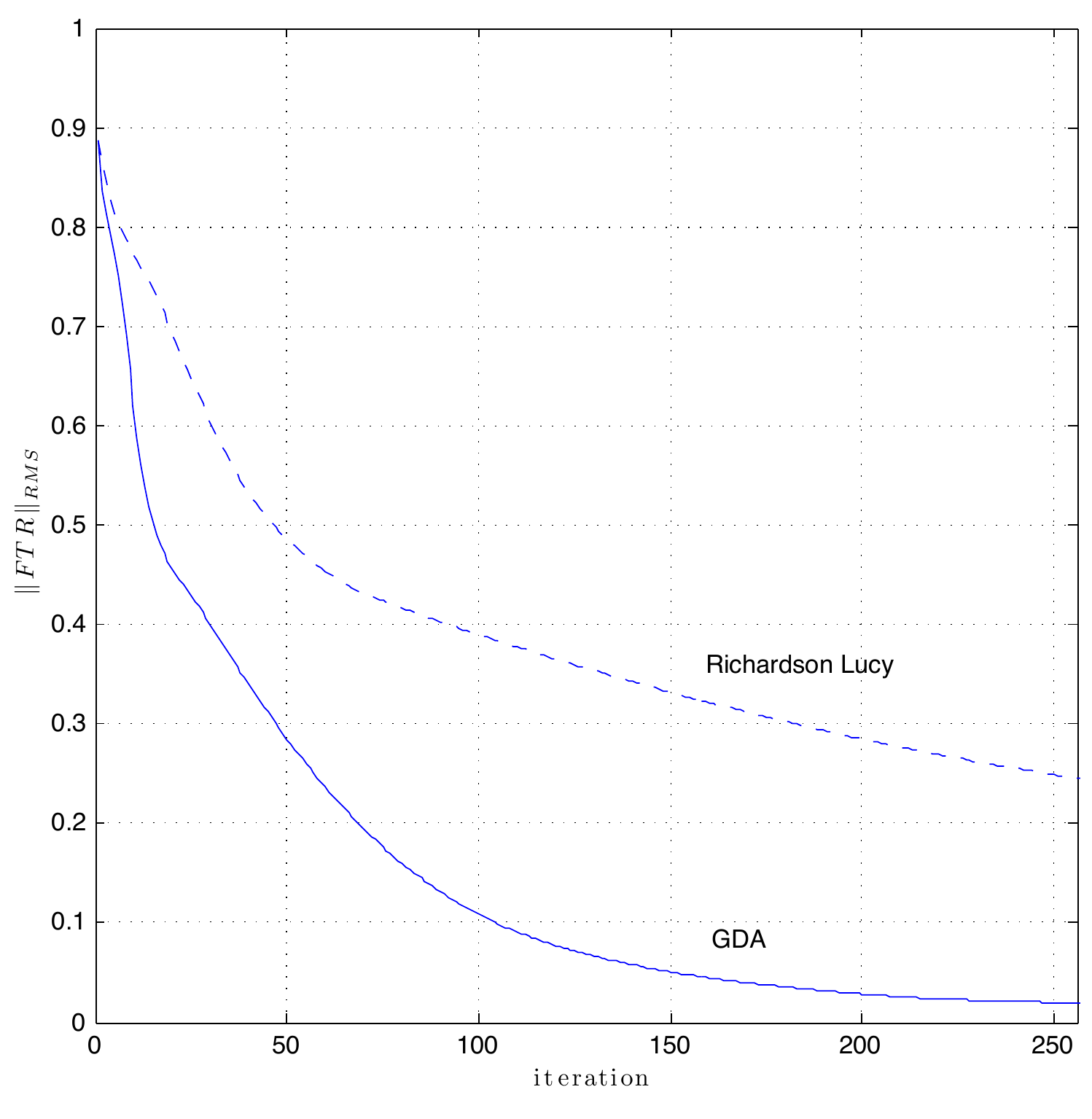}
\par\end{centering}

\caption{(Left) The relative error plotted as a function of iteration for both
the equation (\ref{eq:31}) and Richardson-Lucy algorithms. (Right)
The FTR plotted as a function of iteration for both the equation (\ref{eq:31})
and Richardson-Lucy algorithms, highlights the strong frequency \textquotedblleft reach\textquotedblright{}
of the algorithm compared to RL for this type of image.}

\label{fig:2}
\end{figure}

The effect of the $\omega_{R}$ parameter is evident in the two test
comparisons reported: figure \ref{fig:1} indicates that Richardson-Lucy
recovers the details around Lena's hat relatively poorly, which shows
up as the large difference in FTR plot between the two algorithms
in figure \ref{fig:2}.

\section{Conclusions\label{sec:6}}

The algorithm presented derives the deconvolution operator for images
as a possible Cauchy sequence of weight functions in the space of
solution images. The initial results verify that there is indeed a
measurable advantage to be had in this adopting this approach.

The results here have been deliberately restricted to Cauchy sequences
using symmetric kernels, but in principle the equations can be potentially
extended to different types of kernel through the use of equation
(\ref{eq:13}).

Furthermore, according to lemma \ref{LS1.1}, the Cauchy sequence
demonstrated in equation (\ref{eq:20}) is not necessarily unique.
It can be shown that a more general , exact form of equation (\ref{eq:5})
or (\ref{eq:15}) is:

\begin{equation}
g=\varrho_{0}h+\sum_{n=0}^{\infty}(\mathcal{I}-k*)^{n}(h-h_{0})\label{eq:35}
\end{equation}

which can be used if an operator form for $(\mathcal{I}-k*)^{n}(h-h_{0})$
is available.

Among the advantages associated with techniques that use convergent
sequences like equation (\ref{eq:4}) is that problems like dimensionality,
mathematical representation and the question of continuous or discrete
data spaces, become technical issues. Instead, an appreciation of
the concepts of convergence and accuracy and their impact on the quality
of the image reconstruction become central to the discussion.

This has happened before with the Richardson-Lucy algorithm which
subsequently has found a wide range of application everywhere. Specialising
the image $h$ to a two dimensional series we were able to derive
an improved two-step, non-local version of it for this paper, which
may be used as a simple but promising alternative.

A related fact is that since deconvolution problems occur frequently,
not just in image processing, but are in fact very common to every
area of physics and engineering, a closed form series like equation
(\ref{eq:15}) and more generally (\ref{eq:35}) can potentially be
of use in problems such as viscoelasticity, diffraction and quantum
mechanics to name but a few other areas of potential application.

From a purely theoretical perspective, the technique of finding an
approximate \emph{closed form} Cauchy deconvolution sequence in direct
space mimics the natural one in Fourier space: The Fourier transform
of $g$ is the ratio of the Fourier transform of $h$ and $k$. This
commodity is what in principal makes that approach appealing: the
formal ``rotation'' of the problem into Fourier space \emph{simplifies}
it. The iterative summing of the Cauchy sequence as in equation (\ref{eq:20})
or when applicable (\ref{eq:15}), provides just such an analogue
in direct space once $\varrho$ has been found.

\pagebreak{}

\section{Figures \label{sec:7}}

\section{Bibliography\label{sec:8}}

\end{document}